\documentclass{article}
\usepackage{epsfig,epic,eepic,latexsym, amssymb}

\input xy
\xyoption{all}
\def\pf{PROOF:$\;$}
\def    \QED    {\hfill\hbox{\hskip 4pt
                \vrule width 5pt height 6pt depth 1.5pt}}
\def\epf{\QED\\}

\newcommand{\nit}{\mbox{$\hbox{\it I\hskip -2pt N}$}}

\newcommand{\C}{ {\cal C}}
\newcommand{\V}{ {\cal V} } 
\newcommand{\R}{[0,\infty]}
\newcommand{\TWO}{{\bf 2 }}
\newcommand{\p}{ {\cal P} }
\newcommand{\Q}{ {\cal Q} }

\newcommand{\Pu}{ {\cal P}_1 }
\newcommand{\Pd}{ {\cal P}_{\omega} }
\newcommand{\Pa}{ {\cal P}_{\aleph}}

\newcommand{\FFu} { Fil_1 }

\newcommand{\FFa} { Fil_{\aleph} }
\newcommand{\F}{{\cal F}}
\newcommand{\Fs}{{\cal F}^s}
\newcommand{\Ba}{\Gamma}
\newcommand{\Bas}{\Gamma^s}
\newcommand{\I}{I}
\newcommand{\cI}{{\cal I}}
\newcommand{\Mm}{ M^- }
\newcommand{\Mp}{ M^+ }
\newcommand{\Ml}{ M^l }
\newcommand{\Mr}{ M^r }
\newcommand{\raM}{ \tilde{M} }
\newcommand{\Liminf}{ lim^- }
\newcommand{\Limsup}{ lim^+ }

\newcommand{\VCat}{\V\mbox{-}Cat}
\newcommand{\Phicoc}{\Phi\mbox{-}\mathit{Cocts}}
\newcommand{\adcpo}{\aleph\mbox{-}\mathit{dcpo}}

\newtheorem{theo}{Theorem}[section]
\newtheorem{defi}[theo]{Definition}
\newtheorem{prop}[theo]{Proposition}
\newtheorem{lemm}[theo]{Lemma}
\newtheorem{coro}[theo]{Corollary}
\newtheorem{rema}[theo]{Remark}
\newtheorem{fact}[theo]{}

\title{Flatness, preorders and general metric spaces}
\author{Vincent Schmitt}
\begin{document}
%
%
%
%
\maketitle
\begin{abstract}
We use a generic notion of flatness in the enriched 
context 
to define various completions of metric spaces --
enrichments over $\R$ -- and
preorders -- enrichments over $\TWO$. We characterize 
the weights of colimits commuting in $\R$ with the 
{\em conical terminal object and cotensors}. These 
weights can be interpreted in metric terms
as peculiar filters, the so-called filters
{\em of type $1$}. This generalizes Lawvere's 
correspondence between minimal Cauchy filters and adjoint modules.
We obtain a metric completion
based on the filters of type $1$ as an instance of the
{\em free cocompletion under a class of weights} defined by Kelly.
Another class of flat presheaves is considered
both in the metric and the preorder context. The 
corresponding completion for preorders is the so-called
{\em dcpo completion}. 
\end{abstract}
\section{Introduction}
\label{intro}
Many mathematical objects
have been fruitfully described as enriched categories:
modules, sheaves \cite{Bet85}, \cite{Wal81}, \cite{Wal82}, 
fibrations \cite{BCSW83} and stalks \cite{Str83-2} but 
also, and more simply, metric spaces and preorders 
\cite{Law73}. Amongst the interesting properties inherited 
from category theory, is a general process of (co)completion 
which consists roughly in adding freely to a category,  
colimits of a given kind. 
It is known that these free cocompletions
capture many classical cases of completions,
as for instance the sheafification and the completion \`a 
la Cauchy of metric spaces.
The above completion process for categories have been 
investigated for long and by various authors
(\cite{Ko67}, \cite{Ko95}, \cite{Wo78}, \cite{Tho82}) 
and a general theorem, due to Kelly, may be found in his book
\cite{Kel82}. Kelly's theorem asserts that the free addition
of {\em weighted} colimits yields indeed a universal construction.
Also the important notion of {\em closed class of weights} or 
{\em Betti's covering} occurred in \cite{Bet85} and was studied 
in detail in \cite{AK88} by Albert and Kelly. In a recent 
paper \cite{KS05}, these classes are considered again
and renamed {\em saturated}. A theorem by Albert and Kelly
asserts that the free cocompletion of 
a small category $A$ under a saturated class $\Phi$ of weights
is the full subcategory of presheaves over $A$ generated by 
objects in $\Phi$. From a practical point of view, 
this result may yield meaningful internal descriptions of 
free cocompletions. 
As mentioned in \cite{KS05}, it is the case that 
any class of {\em flat} weights is saturated.
Let us underline an important point: the definition of flatness in 
\cite{KS05} and used here, relies on classes of {\em weights} - 
not on classes of {\em diagrams}, like for instance in \cite{ABLR02}.
A few nice results for this generic notion of flatness 
are recalled in section \ref{flatness}.
The purpose of this paper is to investigate a few 
categorical cocompletions corresponding to various 
notions of flatness in the somewhat simple setting of 
metric spaces -- enrichments over $\R$ --
and preorders -- enrichments over $\TWO$.\\


The results presented in this paper are the following.
We studied the notions of flatness 
based on the following classes of weights:
\begin{itemize}
\item $\Pu$ of weights with domains, either 
the empty category, or the unit category $\I$ (with one point $*$,
and $\I(*,*) = \I$ the unit of $\V$); 
\item $\Pa$ of weights with domains with less
than $\aleph$ objects, for an infinite regular 
cardinal $\aleph$.
\end{itemize}
In the context of metric spaces, or equivalently 
enrichments over the base $\R$, 
$\Pu$-flat presheaves on a small $A$ (seen as a category)
correspond to particular filters, on 
$A$ (seen as a metric space). We named these filters 
{\em the filters of type $1$}.
This generalizes the fact that minimal Cauchy filters 
on $A$ are in one-to-one correspondence 
with the left adjoint modules
$\xymatrix{\I \ar[r]|{\circ} & A}$.
As a consequence of this, one may use
results from (enriched) category
theory to develop a theory of convergence for 
general metric spaces based on the filters of type $1$.
For instance we can derive a completion 
of general metric spaces that we formulate in pure 
topological/metric terms. This completion
may be described as a non-symmetrical version of  
the completion {\em \`a la Cauchy}, the
so-called {\em bi-completion} in \cite{FL82}, 
but it differs from the latter, even for symmetric spaces. 
Secondary completions can also be deduced by considering 
subclasses of filters, namely the
filters of type $\aleph$, which are the metric 
counterparts of the $\Pa$-flat weights.
These latter completions coincide with the 
Cauchy-completion for symmetric spaces. They
also have a ``domain-theoretical'' flavor:   
in the context of preorders and for the case
$\aleph = \omega$, the completion obtained is 
the well-known ``dcpo completion''.\\


The topics of general metric spaces has gained
popularity since Lawvere's papers (\cite{Law73}
and \cite{Law86}), especially
amongst Computer Scientists (\cite{Fla92},
\cite{BvBR98}, \cite{Fla97}, \cite{Wag94}, 
\cite{Rut96}, 
\cite{BvBR98}, \cite{KuSc02}, \cite{Vic} -- recently published 
in \cite{Vic05}).
Surprisingly a few of these works really treat general metric spaces 
explicitly for what they are: enriched categories of a 
particular kind. This is the view in \cite{BvBR98}
where completions are defined
by considering {\em ordinary} colimits in the presheaf 
categories. Also Vickers' recent published 
work \cite{Vic05} concerns actually a peculiar class of flat presheaves,
namely the flat ones with respect to finite conical
weights. It is worth comparing his work with the present 
one. (Both Steve Vickers and the author of the paper
were aware of each others' work.)
We consider in this paper classes of {\em weights} of limits
and colimits and in particular the class of weights
for {\em cotensors}. This class 
is actually crucial to 
characterize the weights corresponding to filters
when the base category is $\R$. 
The approach in this work tries to be 
``as categorical as possible''. The author also hopes that 
it reveals convincing with respect to non-categorical 
issues. For instance it helps to develop a theory of convergence
for non-symmetric space. It may also precise the essence of
the ``similarity'' between directed complete partial 
orders and complete metric spaces: both being after all, 
the algebras of very similar {\em KZ-doctrines}.\\

The paper is organized as follows.
Section \ref{flatness} briefly recalls 
some material about enriched category theory. It
aims at precising the terminology,
and presenting most of the categorical concepts and 
results used in the rest of the paper.
Next sections are devoted
to applications, namely
general metric spaces (the base $\V$ is $\R$) in section \ref{gms} 
and preorders ($\V$ is $\TWO$) in section \ref{preorder}.\\

I thankfully acknowledge Claudio Pisani, Pierre Ageron and Isar Stubbe
for their helpful comments.\\

\section{Background in Enriched Category Theory}
\label{flatness}
For the background
knowledge about enriched categories, we refer the reader to 
Kelly's book \cite{Kel82}, augmented by Betti's \cite{Bet85}, 
the Albert-Kelly article \cite{AK88} and the recent 
paper \cite{KS05}.\\

We consider in this section, a symmetric monoidal
complete and cocomplete closed $\V$. 
As usual, weights are  
presheaves over {\em small} categories.
Remember that given a family $\Phi$ of weights, 
and a small category $A$, the 
{\em closure $\Phi(A)$ of $A$ under $\Phi$-colimits in $[A^{op},\V]$},  
is defined as the smallest full replete subcategory 
of the presheaf category $[A^{op}, \V]$, containing the representables and 
closed under the formation of $\Phi$-colimits in $[A^{op},\V]$.
The latter means that for any weights $F: K^{op} \rightarrow \V \in \Phi$,
and any functor $G: K \rightarrow [A^{op},\V]$ taking its values 
in $\Phi(A)$, $F * G$ is again in $\Phi(A)$. Remember also that for 
a small $A$, $\Phi(A)$ is generally not small. 
\begin{theo}(\cite{Kel82} Theorem 5.35)
\label{freecoc} For any family of weights $\Phi$,
and for any small $A$, the closure $\Phi(A)$ of $A$
in $[A^{op},\V]$ under $\Phi$-colimits 
constitutes the {\em free ${\Phi}$-cocompletion} of $A$.   
Precisely this means that:
\begin{itemize}
\item 
$\Phi(A)$ is $\Phi$-cocomplete;
\item For any possibly large $\Phi$-cocomplete 
category $B$, one has an equivalence
\begin{center} $(*)$
$Lan_K: \VCat(A,B) \cong \Phicoc(\Phi(A), B)$
\end{center}
\end{itemize} 
where:
\begin{itemize}
\item $K$ is the full and faithful inclusion 
$A \rightarrow \Phi(A)$ sending any $a \in A$ to $A(-,a)$;
\item $\Phicoc(\Phi(A), B)$ is the full subcategory
of $\VCat(\Phi(A),B)$ generated by the
$\Phi$-cocontinuous functors;  
\item $Lan_K$ stands for the {\em left Kan extension functor} and
has inverse the functor ``$- \circ K$'' given by the
composition with $K$.
\end{itemize}
\end{theo}
Actually the equivalence $(*)$ above lifts to a
$\V$-equivalence since the $A$ considered is small.

\begin{defi}
Given a class $\Phi$ of weights, its {\em saturation}
is the largest class $\Phi^*$ such that
every $\Phi$-complete category is also $\Phi^*$-complete, and
every $\Phi$-continuous functor is also $\Phi^*$-continuous.
\end{defi}
Note that in the definition above, one can substitute 
``$\Phi$-cocomplete'' to ``$\Phi$-complete'' and 
``$\Phi$-cocontinuous'' to ``$\Phi$-continuous''. 
A class $\Phi$ of weights 
is called {\em saturated} if and only if $\Phi = \Phi^*$.
For any small $\V$-category $A$, and any family $\Phi$ of weights,
we write
$\Phi[A]$ for the full subcategory 
of $[A^{op},\V]$ with objects in $\Phi$.
Note that by Yoneda, any weight $F$ occurs as a colimit
weighted by $F$ of representables: $F \cong F*Y$,
so that for any class $\Phi$ of weights and 
any small category $A$, certainly $\Phi[A] \subset \Phi(A)$.\\

We can now formulate the Albert-Kelly theorem \cite{AK88}.
\begin{theo}\label{AK} 
For any class $\Phi$ of weights, a 
weight $\phi: A^{op} \rightarrow \V$ lies in 
$\Phi^*$ if and only if 
$\phi \in \Phi(A)$. Which is to say that $\Phi^*[A] = \Phi(A)$
for all class $\Phi$ and all small $A$.
\end{theo} 
\begin{coro} 
A class $\Phi$ is saturated if and only if 
$\Phi[A] = \Phi(A)$ for all small category $A$.
\end{coro} 
It is our view that the latter characterization of saturated 
classes can reveal useful to obtain internal descriptions of
free cocompletions. As shown in \cite{KS05},
an important example of saturated classes
is given by classes of {\em flat} weights that we introduce 
now.\\

Given a class $\Phi$ of weights, a weight $F: A^{op} \rightarrow \V$ 
is said {\em $\Phi$-flat} when its left Kan extension along $Y$,
$-*F : [A,\V] \rightarrow \V$ preserves all $\Phi$-limits.
This is also equivalent to say that the 
``weighting by $F$'' functor $F*- : [A,\V] \rightarrow \V$
-- that is isomorphic to $-*F$ -- is
$\Phi$-continuous. As precisely stated in \cite{KS05},
$F$ is $\Phi$-flat when {\em $\Phi$-limits commute with $F$-colimits 
in the base $\V$}. Letting $\Phi^+$ denote the family of
$\Phi$-flat weights, it happens that the class $\Phi^+$ is
saturated. Thus by the Albert-Kelly theorem one has:
\begin{theo}
\label{fcomp}
For any family $\Phi$ of weights and 
any small category $A$, $\Phi^+[A]$ is the free 
$\Phi^+$-cocompletion $\Phi^+(A)$ of $A$.
\end{theo}

We wish to mention here another nice consequence of the
fact that classes of the form $\Phi^+$ are saturated. 
This is obtained by combining Theorem \ref{fcomp} above 
with the characterization (due to Kelly) of the embeddings of the 
form $A \rightarrow \Phi(A)$ (see again \cite{KS05}, Propositions
4.2 and 4.3).
For a class $\Phi$ and a category $A$, we write 
$A_{\Phi}$ for its full subcategory defined by the
$a$ such that $A(a,-): A \rightarrow \V$ preserves 
$\Phi$-colimits. 
\begin{theo}
For any class $\Phi$ of weights, 
if $\Psi$ stands for $\Phi^+$, then
for a category $A$, the following are equivalent:\\  
- $A$ is $\Psi$-cocomplete and has a full small category 
$B \subset A_{\Psi}$ such that any  $a \in A$ is a 
$\Psi$-colimit of a diagram in $B$;\\
- there exists a small $B$ such that $A$ is equivalent to
the full subcategory of $[B^{op},\V]$ defined by $\Phi$-flat weights.
\end{theo}
The previous theorem can be seen as a generalization
of a well-known result regarding accessible categories.
The work \cite{KS05} contains a few
examples of notions of flatness. Let us recall in particular that
for the empty class $0$ of weights and the class 
$\p$ of all weights, one has $0^+ = \p$ whereas
$\p^+$ is the class denoted $\Q$ of {\em small projectives},
which are also the weights $A^{op} \rightarrow \V$
defining left adjoint modules $\xymatrix{I \ar[r]|{\circ} & A}$,
and also the {\em weights of absolute colimits} (\cite{Str83-1}). 
So actually for any small category $A$, $\p^+(A) = \Q(A)$ 
is the Lawvere-Cauchy-completion $A$.\\

In this work, we shall investigate a few classes of 
colimits commuting with certain given limits in $\R$
and in $\TWO$. 
Precisely we shall study the notions of flatness
associated to the
the following classes of weights (that we morally think 
as weights of limits): 
\begin{itemize}
\item $\Pu$: the class of weights of the form 
$F:K \rightarrow \V$ where $K$ is 
the empty $\V$-category or $K = \I$ (the unit category 
with one point $*$ and $\I(*,*) = I$); 
\item $\Pa$, for any infite regular cardinal
$\aleph$:
the class of weights $F:K \rightarrow \V$
such that  and 
$\sharp Obj(K) < \aleph$.
\end{itemize}

We need now to consider ``special'' weighted limits, namely 
{\em cotensors} and {\em conical limits}.  
We skip their definitions, referring the reader to the section 3
of Kelly's book \cite{Kel82}. We recall anyhow a few facts 
regarding {\em conical limits} to 
avoid a too common confusion with the mere ordinary limits.
Conical limits (and dually colimits) are defined as 
representatives of particular presheaves and 
their defining isomorphisms are {\em $\V$-natural}. Though any 
conical limit has an underlying ordinary limit, the converse does 
not hold. There is an alternative characterization of canonical 
limits: they are the ordinary limits which are preserved by 
ordinary functors underlying representables (see \cite{Kel82} p.95).
This point is crucial for instance in the proof of Theorem 
\cite{Kel82} (3.73) and thus in the proof below of Proposition \ref{ect1}.\\ 

According to Kelly's terminology, $\Pu$ is the class
of weights for the conical terminal object and cotensors.
One has the following inclusions of classes of weights, 
where $\aleph$ stands for any infinite regular cardinal: 
\begin{fact}
$\Pu \subset \Pd \subset \Pa \subset \bigcup_{\aleph} \Pa = \p$
\end{fact}
and thus
\begin{fact}\label{hierarchie}
$\Pu^+ \supset \Pd^+ \supset \Pa^+ \supset 
\bigcap_{\aleph} \Pa^+ = \Q$.   
\end{fact}

For any infinite regular cardinal $\aleph$, 
we call a limit with weights in $\Pa$,
an {\em $\aleph$-limit}. 
Thus a {\em conical $\aleph$-limit} 
is just a conical limit whose underlying ordinary limit has an 
indexing diagram with less than $\aleph$ objects.
A minor adaptation of the proof of theorem \cite{Kel82} (3.73) 
like in \cite{Kel82-2} (4.3), shows the following. 
\begin{prop} 
\label{ect1}
For any infinite regular cardinal $\aleph$, a
$\V$-category $A$ is $\Pa$-complete if and only if it 
has cotensors and all conical $\aleph$-limits. 
Given a $\Pa$-complete $A$, a $\V$-functor $P:A \rightarrow B$ 
is $\Pa$-continuous if only if it preserves conical $\aleph$-limits 
and cotensors.
\end{prop}
\pf (Sketch of) It suffices to reuse the argument
developed in the sketch of proof of \cite{Kel82-2} (4.3). 
Remark that if $A$ has all conical $\aleph$-limits and cotensors, 
then for any $F: K \rightarrow \V$ and $G: K \rightarrow A$ 
with $\sharp Obj(K) < \aleph$, the weighted limit $\{F,G\}$  
may be computed as the equalizer in $A_0$ of the canonical 
pair 
$$\xymatrix{ \prod_{k \in K} Fk \pitchfork Gk \ar@<1ex>[r]
\ar@<-1ex>[r] &
\prod_{k,k' \in K} K(k,k') \pitchfork (Fk
\pitchfork Gk')}.$$ 
Actually all the ordinary limits involved 
in this equalizer, i.e. the two products and the equalizer itself, 
have indexing diagrams with less than $\aleph$ objects
and thus are conical.  
Also revisiting the sketched proof of theorem (3.73) in \cite{Kel82},
one gets that for any functor $H: A \rightarrow B$ preserving
conical $\aleph$-limits and cotensors, $H$ will preserve  
an equalizer as above which image in $B$ is then the limit 
$\{F, HG\}$.
\epf

We finish this quick tour of category theory
by mentioning miscellaneous results that will 
serve later.
The following proposition already occurs in \cite{KS05}
but under a more elaborate form.
\begin{prop}
\label{Lprpf} For any saturated class $\Psi$ of weights,
any weight $F: A \rightarrow \V$ in $\Psi$,
and any functor $G: A \rightarrow B$ with $B$ small, the 
left Kan extension of $F$ along $G$ is again in $\Psi$.
\end{prop}
\pf Since $F \cong F*Y$, $F$ is a $\Psi$-colimit of representables.
The image by $Lan_G: [A,\V] \rightarrow [B,\V]$ of any representable 
is again representable  
(For any $a \in A$, $Lan_G(A(a,-))(b) \cong_b B(G-,b) * A(a,-) 
\cong_b B(Ga,b)$).  
Also the left Kan extension functor 
$Lan_G : [A, \V] \rightarrow [B, \V]$ is cocontinuous since
it is a left adjoint, so 
$Lan_G(F)$ is also a $\Psi$-colimit of representables
hence it belongs to $\Psi$ since $\Psi$ is saturated.
\epf

Another point is that when the base $\V$ is 
small -- and this is the case for our applications
with $\V= \TWO$ and $\V = \R$ --
then $\V$ being small-complete is necessarily a preorder 
(see for instance \cite{Bor94-1} prop. 2.7.1 p.59). 
So in the case of a small $\V$, for any small
$\V$-category $A$, the presheaf category $[A,\V]$ remains
small and so does $\Phi(A)$ for any family $\Phi$ of weights.\\

Recall also that if a 
small category $A$ is $\Phi$-cocomplete then
it is a retract of $\Phi(A)$ (i.e. the inclusion 
$A \rightarrow \Phi(A)$ is a split monic) 
but it is generally NOT
isomorphic to $A$. This situation has been
studied (see for instance \cite{Bet85} Theorem p.175, 
or \cite{AK88} Proposition 4.5)
and it happens that a category $A$ is
$\Phi$-cocomplete if and only if the inclusion
$A \rightarrow \Phi(A)$ is a left adjoint.
Nevertheless for some $\Phi$ and some $\Phi$-cocomplete 
$A$, it happens sometimes that $\Phi(A) \simeq A$.
This is the case for all $\Phi$-cocomplete $A$ when 
$\Phi$ is the class $\Q$ of small projective weights. 
As we shall see now, this 
also happens for the base $A = \V$ when 
$\Phi$ is class $\C$ of weights of 
cotensors -- i.e. those weights of the form $\I \rightarrow \V$.
\begin{prop}\label{flacot}
For any monoidal closed complete and cocomplete $\V$,
$\V \simeq \C^+(\V)$.
\end{prop}
\pf
Let us establish first the following result.
\begin{lemm}
Let $M$ and $N$ be two accessible presheaves $\V^{op} \rightarrow \V$,
such that $N$ is $\C$-flat, then $\{ M, N \} = [M*1,N*1]$.
\end{lemm}
\pf
For any $v \in \V$, one has
$[M*1,v] \cong \{ M, [-,v] \}$ so that 
$[M*1,N*1] = \{ M, [-,N*1] \}$, the right hand
side here is $\{ M_{-}, N_? * [-,?] \} = \{ M, N \}$
since $N$ is $\C$-flat and colimits are pointwise
in functor categories.  
\epf

As a consequence of the previous lemma,
the restriction
of $-*1: [\V^{op},\V] \rightarrow \V$
to $\C^+(\V)$ is fully faithful.
Since this functor is also essentially surjective
on objects, 
it is part of an equivalence 
of categories,  
with inverse the inclusion $\V \rightarrow [\V^{op},\V]$.
\epf

\section{General metric spaces -- The case $\V = \R$.}
\label{gms}
This section applies the categorical results presented
in the previous section in the context of general metric spaces.
Let us start by recalling briefly a few results that 
are from \cite{Law73} or belong to folklore.
$\R$ stands for the monoidal closed category with:
\begin{itemize}
\item objects: nul or positive reals and $+\infty$;
\item arrows: the reverse ordering, $x \rightarrow y$ if and 
only if $x \geq y$;
\item tensor: the addition (with
$+\infty + x = x + +\infty =  +\infty$); 
\item unit: 0.
\end{itemize}
For any pair $x,y$ of objects in $\R$, the exponential object 
$[x,y]$ is $max\{ y - x, 0\}$.\\

A {\em small} {\em $\R$-category} $A$ corresponds to a {\em general metric space}.
That is a set of {\em objects} or {\em elements}, $Obj(A)$ 
(most of the time just denoted by $A$) together with a 
map $A(-,-):  Obj(A) \times  Obj(A) \rightarrow \R$,
called {\em pseudo-distance},
that satisfies:
\begin{itemize}
\item for all $x,y,z \in A$, $A(y,z) + A(x,y) \geq A(x,z)$; 
\item for all $x \in A$, $0 \geq A(x,x)$.
\end{itemize}
A $\R$-functor $F: A \rightarrow B$
corresponds to a {\em non-expansive} map $F: Obj(A) \rightarrow Obj(B)$, i.e.
for all $x,y \in A$, $A(x,y) \geq B(F(x), F(y))$. 
A $\R$-natural transformation $F \Rightarrow G: A \rightarrow B$
corresponds to the fact that
for all $x \in  A$, $0 \geq B(F(x), G(x))$.
A $\R$-module $M: \xymatrix{I \ar[r]|{\circ} & A}$ -- or left module
on $A$ -- is a map 
$Obj(A) \rightarrow \R$ such that for all $x,y \in A$,
$M(y) + A(x,y) \geq M(x)$.
Dually a $\R$-module $N: \xymatrix{A \ar[r]|{\circ} &  I}$ -- or right module 
on $A$ -- is a map 
$Obj(A) \rightarrow \R$ such that for all $x,y \in A$,
$A(x,y) + N(x) \geq N(y)$.
For any general metric spaces $A$ and $B$, the set
of non expansive maps from $A$ to $B$ becomes 
a general metric space $[A,B]$ with distance 
$[A,B](f,g) = \bigvee_{x \in A} B(f(x),g(x))$.
In particular the presheaf category $[A^{op},\R]$ has homsets
given by $[A^{op},\R](M,N) = \bigvee_{x \in A}[M(x),N(x)]$.
Its underlying category is a partial order
with arrows given by the pointwise reverse ordering:
$M \Rightarrow N$ if and only if $\forall x \in A$,
$M(x) \geq N(x)$.
The composition of left and right modules is as follows.
Given $\xymatrix{ I \ar[r]|{\circ}^M & A \ar[r]|{\circ}^N & I}$, the
composite $N * M$ is $\bigwedge_{x \in A} M(x) + N(x)$.
For such $M$ and $N$,
$M$ is left adjoint to $N$
if and only if:
\begin{itemize}
\item $(1)$ $0 \geq N * M$; 
\item $(2)$ for all $x,y \in A$, $N(y) + M(x)  \geq A(x,y)$.
\end{itemize}

For the rest of this section, $A$ denotes an arbitrary
general metric space that we freely see as a small 
category.\\

Lawvere observed that Cauchy sequences on the space
$A$ correspond to (left) adjoint modules $\xymatrix{I \ar[r]|{\circ} & A}$.
Actually there is a 
bijection between left adjoint modules on the small $\R$-category 
$A$ and minimal Cauchy filters on the space $A$.
From this observation mainly, one gets that the full 
subcategory of $[A^{op},\R]$ generated by left adjoint 
modules $\xymatrix{ \I \ar[r]|{\circ} & A}$, 
corresponds to the completion \`a la Cauchy of the space
$A$ if $A$ is symmetric
or in general to its bi-completion
(see for instance \cite{FL82} and \cite{Fla92}
or \cite{Sch03}).
\begin{defi}\label{defCauchy}
A filter $\F$ on $A$ is {\em Cauchy} if 
and only if for any $\epsilon > 0$, there exists an $f \in \F$
such that for any elements
$x, y$ of $f$, $A(x,y) \leq \epsilon$ or equivalently when: 
$$\bigwedge_{f \in \F} \bigvee_{x,y \in f} A(x, y)= 0.$$
\end{defi}
\begin{defi}
For any left adjoint module $M$ on $A$, with right adjoint $\raM$,
$\Bas(M)$ stands for the subset of 
the powerset $\wp(A)$ of $A$:
$\{ \Bas(M)(\epsilon) \mid \epsilon \in ]0, + \infty] \}$,
where $\Bas(M)(\epsilon)$ denotes the set $\{ x \in A \mid 
M(x) + \raM(x) \leq \epsilon\}$.  
\end{defi}
For any left adjoint module $M$ on $A$, $\Bas(M)$ is a Cauchy 
basis. The filter that it generates, that we denote $\Fs(M)$, is a 
minimal Cauchy filter. The map $M \mapsto \Fs(M)$ defines a
bijection between left adjoint 
modules $\xymatrix{ I \ar[r]|{\circ} & A}$
and minimal Cauchy filters on $A$.
One may check the following points (proved for instance in \cite{Sch03}).
To any Cauchy filter $\F$, one may associate 
a left adjoint module $\Ml(\F)$ defined by
$$x \mapsto \bigwedge_{f \in \F} \bigvee_{y \in f} A(x,y) 
         = \bigvee_{f \in \F} \bigwedge_{y \in f} A(x,y).$$
$\Ml(\F)$ has right adjoint $\Mr(\F)$ given by the map
$$x \mapsto \bigwedge_{f \in \F} \bigvee_{y \in f} A(y,x)= 
\bigvee_{f \in \F} \bigwedge_{y \in f} A(y,x).$$
For any left adjoint module $M$ on $A$, 
$\Ml(\Fs(M)) = M$ and
for any Cauchy filter $\F$ on $A$, 
$\Fs(\Ml(\F))$ is the only
minimal Cauchy filter contained in $\F$.


\subsection{Modules and Filters}
Since the (Lawvere-)Cauchy-completion of $A$
is up to equivalence the full subcategory
of flat presheaves $\p^+(A)$, 
we wondered whether the previous correspondence between 
left adjoint modules and Cauchy filters could be extended to 
larger classes of $\Phi$-flat modules and filters. 
We are going to show in \ref{mainres} that this is the case for 
$\Phi$ = $\Pu$ 
as we exhibit a $\R$-category structure 
$\FFu(A)$ on the set of the so-called {\em filters of type $1$} 
on $A$, which is equivalent to the category $\Pu^+(A)$.\\ 

We need to recall a few technical points before
giving simple characterizations of the 
$\Pu$ and $\Pa$-flat presheaves.
For the assertions \ref{cotpoint}, \ref{precon} and 
\ref{colimpre} below, $\V$ denotes a complete and cocomplete
monoidal closed category.
Remember that cotensors are defined pointwise in functor
categories. In particular: 
\begin{fact}
\label{cotpoint} For any category $C$,
the presheaf $\V$-category $[C,\V]$ is cotensored and for any 
presheaf $N$,
$v \pitchfork N$ is the composite 
$\xymatrix{ C \ar[r]^{N} & \V \ar[r]^{[v,-]} & \V }$.
\end{fact}
Also for functors between {\em cocomplete} categories,
the preservation of conical colimits amounts
to the preservation of ordinary colimits. Precisely
one may check:
\begin{fact}
\label{precon}
Given a $\V$-functor $T:C \rightarrow D$ 
with underlying ordinary functor
$T_0: C_0 \rightarrow D_0$ 
and an ordinary functor $P: J \rightarrow C_0$ with $J$ 
small, if the conical limits of $P$ and of $T_0P$ exist
and $T_0$ preserves the ordinary limit of $P$,
then $T$ preserves the conical limit of $P$.
\end{fact}
Eventually the preservation of limits/colimits is simple
in the case $\V = \R$, according to the following observation.
\begin{fact}
\label{colimpre}
If the base category $\V_0$ is a preorder, then
given a weight $F: B^{op} \rightarrow \V$,
a functor $G: B \rightarrow C$ such that $F*G$ exists
and a functor $H: C \rightarrow D$, $H$ preserves $F*G$ if and
only if $F*(GH)$ exists and $H(F*G) \cong F*(GH)$.
\end{fact}

According to the three previous points one gets the following
results.
\begin{fact}
Let $M: \xymatrix{I \ar[r]|{\circ} & A}$ be a left module.
\begin{itemize}
\item $-*M: [A,\R] \rightarrow \R$ preserves 
the unique conical limit with weight with empty
domain if and only if the underlying ordinary functor 
preserves the terminal object
i.e. $0 * M = 0$ if and only if 
$$(1)\;\;\bigwedge_{x \in A} M(x) = 0.$$  
\item For an infinite regular cardinal $\aleph$,
$-*M$ preserves the conical $\aleph$-limits if and only if\\ 
$(2)$ for any family of right modules 
$N_i: \xymatrix{A \ar[r]|{\circ}  & I}$, 
where $i \in I$ and $\sharp I < \aleph$,
$$\bigwedge_{x \in A}( M(x) + \bigvee_{i \in I} N_i(x) )
= \bigvee_{i \in I}( \bigwedge_{x \in A} M(x) + N_i(x)  );$$
\item $-*M$ preserves cotensors if and only if\\
$(3)$ for any $v \in \R$ and 
any right module $N: \xymatrix{A \ar[r]|{\circ} & I}$ , 
$$\bigwedge_{x \in A}( M(x) + [v,N(x)] ) = [v, \bigwedge_{x \in A}( M(x) + N(x) )].$$
\end{itemize}
\end{fact}
So $\Pu$-flat modules are the modules satisfying $(1)$ and $(3)$ above,
whereas for an infinite regular cardinal $\aleph$,
the $\Pa$-flat modules are those satisfying $(2)$ and $(3)$.\\  

We introduce now the filters that will
occur as the metric counterparts of the $\Pu$-flat
modules.
\begin{defi}
Given a filter $\F$ on $A$ and a map $f: Obj(A) \rightarrow Obj(\R)$,
$\Limsup_{x \in \F} f(x)$, or simply $\Limsup_{\F} f$, denotes 
$\bigwedge_{f \in \F} \bigvee_{x \in f} f(x)$ and 
 $\Liminf_{x \in \F} f(x)$, or $\Liminf_{\F} f$, denotes
$\bigvee_{f \in \F} \bigwedge_{x \in f} f(x)$.
\end{defi}
\begin{defi}\label{deftype1}
A filter $\F$ on $A$ has
{\em type 1} if and only if 
$$\Limsup_{x \in \F} \Liminf_{y \in \F} A(x,y) = 0.$$
\end{defi}
Now compare the previous definition and \ref{defCauchy}.
Definition \ref{deftype1} is a generalization to 
non-symmetric 
spaces of the fact that the diameter of the elements of 
the filter $\F$ may be chosen arbitrary small. 
\begin{rema}
Any Cauchy filter has type $1$.
\end{rema}

From the correspondence between Cauchy filters and left adjoint modules,
we know two operators associating filters to modules.
\begin{defi}
Given any filter $\F$ on $A$, we define the following $\R$-valued maps
on objects of $A$: 
\begin{tabbing}
\hspace{1cm}\=$\Mm(\F): x \mapsto \Liminf_{\F} A(x,-) 
= \bigvee_{f \in \F} \bigwedge_{y \in f} A(x,y)$,\\
\>$\Mp(\F): x \mapsto  \Limsup_{\F} A(x,-) = \bigwedge_{f \in \F} \bigvee_{y \in f} A(x,y)$.
\end{tabbing}
\end{defi}
For any filter $\F$ on $A$, one has $\Mm(\F) \leq \Mp(\F)$ where 
the order is pointwise, and if $\F$
is Cauchy then $\Mm(\F) = \Mp(\F)$. 
\begin{rema}
A filter $\F$ on $A$ has
{\em type 1} if and only if $$\Limsup_{\F} \Mm(\F) = 0.$$
\end{rema}

\begin{fact}
Given any filter $\F$ on $A$,
the map $x \mapsto \Mm(\F)(x)$ defines a module 
$\xymatrix{I \ar[r]|{\circ} & A}$.
\end{fact}
\pf
For all $x,y \in A$, 
\begin{tabbing}
$\Mm(\F)(x) + A(y,x)$ \=$=$ \=$( \bigvee_{f \in \F} \bigwedge_{z \in f}  A(x,z) ) + A(y,x)$\\
\>$\geq$ \>$\bigvee_{f \in \F} ( (\bigwedge_{z \in f} A(x,z)) + A(y,x)   )$\\
\>$=$\>$\bigvee_{f \in \F} \bigwedge_{z \in f} (A(x,z) + A(y,x))$\\
\>$\geq$\>$\bigvee_{f \in \F} \bigwedge_{z \in f} A(y,z)$\\
\>$=$ $\Mm(\F)(y)$.
\end{tabbing}
\epf
Conversely assigning a filter to a module should be a simple matter. 
\begin{defi}
For any module $M$, $\Ba(M)$ denotes the subset  
$\{ \Ba(M)(\epsilon) \mid \epsilon \in ]0, + \infty] \}$
of $\wp(A)$,
where $\Ba(M)(\epsilon)$ is the set $\{ x \in A \mid 
M(x) \leq \epsilon\}$. Also $\F(M)$ denotes the upper 
closure of $\Ba(M)$ in $\wp\wp(A)$ ordered by inclusion.
\end{defi}



We can now state our main result:
\begin{theo}\label{mainres}
The set of filters of type $1$ on $A$ 
may be given a general metric space structure $\FFu(A)$
that is equivalent to $\Pu^+(A)$.  The ``distance'' on 
$\FFu(A)$ is defined by the map 
$$(\F_1,\F_2) \mapsto \Limsup_{x \in \F_1} \Liminf_{y \in \F_2} A(x,y).$$ 
The functors of this equivalence are defined by the 
maps on objects $M \mapsto \F(M)$ and $\F \mapsto \Mm(\F)$.
\end{theo}
We prove now a succession of results that constitutes 
the proof of Theorem \ref{mainres}.

\begin{prop}
\label{fac22} For any module $N : \xymatrix{ A \ar[r]|{\circ} &  I}$,
and any filter $\F$ on $A$,
$$N * \Mm(\F) \geq \Liminf_{\F} N.$$
Moreover if $\F$ is of type 1 then
the previous inequality becomes an equality.
\end{prop}
\pf
For any module $N$ and any filter $\F$ as above,
\begin{tabbing}
\hspace{1cm}$N*\Mm(\F)$\=$=$\=$\bigwedge_{x \in A} ( \Mm(\F)(x) + N(x) )$\\
\>$=$ \>$\bigwedge_{x \in A}( ( \bigvee_{f \in \F} \bigwedge_{y \in f} A(x,y)  )  + N(x)  )$\\
\>$\geq$ \>$\bigwedge_{x \in A} \bigvee_{f \in \F} ( (\bigwedge_{y \in f} A(x,y))    + N(x))$\\
\>$=$ \>$\bigwedge_{x \in A} \bigvee_{f \in \F} \bigwedge_{y \in f} ( A(x,y)    + N(x))$\\
\>$\geq$ \>$\bigvee_{f \in \F} ( \bigwedge_{y \in f} N(y) )$.
\end{tabbing}
Let us suppose moreover that $\F$ has type 1.
Let $\epsilon > 0$. One may choose $f_{\epsilon} \in \F$ such that
when $x \in f_{\epsilon}$, $\Mm(\F)(x) \leq \epsilon$.
Thus \begin{tabbing}
\hspace{1cm}$N * \Mm(\F)$\=$=$\=$\bigwedge_{x \in A} ( \Mm(\F)(x) + N(x) )$\\ 
\>$\leq$ \>$\Mm(\F)(x) + N(x)$, for any $x \in f_{\epsilon}$\\
\>$\leq$ \>$\epsilon + N(x)$, for any $x \in f_{\epsilon}$.
\end{tabbing}
Thus 
\begin{tabbing}
\hspace{1cm}$N * \Mm(\F)$\=$\leq$\=
$\bigwedge_{x \in f_{\epsilon}} ( \epsilon + N(x) )$\\
\>$=$ 
\>$\epsilon + \bigwedge_{x \in f_{\epsilon}} N(x)$\\
\>$\leq$ 
\>$\epsilon + \bigvee_{f \in \F} \bigwedge_{x \in f} N(x)$.     
\end{tabbing}
\epf

\begin{fact}
\label{cot}
If $\F$ is a filter of type 1 on $A$ then 
$ - * \Mm(\F)$ preserves cotensors.
\end{fact}
The proof of this fact relies on a very peculiar
property of the base $\R$, namely that {\em cotensors commute 
in $\R$ with conical colimits of non-empty diagrams}. 
This is to say:
\begin{fact}
\label{facR}
For any $v$ in $\R$ and any non empty family $(a_i)_{i \in I}$ in $\R$,
$$[v,\bigwedge_{i \in I} a_i] = \bigwedge_{i \in I} [v,a_i].$$
\end{fact}
We just need to prove 
$[v,\bigwedge_{i \in I} a_i] \geq \bigwedge_{i \in I} [v,a_i]$.
Let us fix $\epsilon > 0$.  
Since $I$ is not empty, there exists $j \in I$ such that
$\epsilon + \bigwedge_{i \in I} a_i \geq a_j$. Also
$[v, \bigwedge_{i \in I} a_i] \geq [v, \bigwedge_{i \in I} a_i]$,
so $v + [v, \bigwedge_{i \in I} a_i] \geq \bigwedge_{i \in I} a_i$.
For a $j$ as above,
$\epsilon + v + [v, \bigwedge_{i \in I} a_i] \geq a_j$
and $\epsilon + [v, \bigwedge_{i \in I} a_i] \geq [v, a_j] 
\geq \bigwedge_{i \in I}[v,a_i]$.  
\epf

We can now prove \ref{cot}.\\
\pf Given $v \in \R$ and $N: \xymatrix{ A \ar[r]|{\circ} & I}$
we have to show $(v \pitchfork N) * \Mm(\F) =  [v, N*\Mm(\F)]$. 
According to \ref{fac22}, 
$(v \pitchfork N) * \Mm(\F)$ $=$ 
$\bigvee_{f \in \F} \bigwedge_{x \in f} [v,N(x)]$
and $[v, N * \Mm(\F) ]$ $=$ 
$[v,  \bigvee_{f \in \F} \bigwedge_{x \in f} N(x)]$
$=$ 
$\bigvee_{f \in \F} [v, \bigwedge_{x \in f} N(x)]$.
Since all the $f \in \F$ are non empty, the result follows then 
from \ref{facR}.
\epf

\begin{fact}
\label{wfpt}
If $\F$ is a filter of type 1 on $A$
then $-* \Mm(\F)$ preserves the terminal object. 
\end{fact}
\pf
We have to show that $\bigwedge_{x \in A} \Mm(\F)(x) = 0$.
For any $\epsilon > 0$, since $\Limsup_{\F}  \Mm(\F) = 0$,
one may find an $f \in \F$ such that for any $x \in f$,
$\Mm(\F)(x) \leq \epsilon$. Since that $f$ is not empty,
$\bigwedge_{x \in A} \Mm(\F)(x) \leq \epsilon$.
\epf

According to \ref{cot} and \ref{wfpt}, for any filter
$\F$ of type $1$, the module $\Mm(\F)$ is $\Pu$-flat.
Therefore we shall define a general metric structure
$\FFu(A)$ on the set of filters of type 1 by letting
$$\FFu(A)(\F_1,\F_2) = [A^{op},\R](\Mm(\F_1),\Mm(\F_2)).$$
(We shall simplify this distance later on.)
A consequence of this definition is that the 
map $\F \mapsto \Mm(\F)$ defines 
a fully faithful functor $\FFu(A) \rightarrow \Pu^+(A)$.

\begin{fact}
\label{unit}
A filter $\F$ on $A$ has type $1$ if and only if
$\F \supset  \F \circ \Mm ( \F )$. 
\end{fact}
\pf One has the successive equivalences.
\begin{tabbing}
\hspace{0.5cm}\=$\F$ has type 1\\ 
\hspace{1cm}if and only if\\
\>$\bigwedge_{f \in \F} \bigvee_{x \in f} \Mm(\F)(x) = 0$\\
\hspace{1cm}if and only if\\
\>for all $\epsilon > 0$, 
there exists $f \in \F$ such that 
for all $x \in f$, $\Mm(\F)(x) \leq \epsilon$,\\
\hspace{1cm}if and only if\\
\>for all $\epsilon > 0$, 
there exists $f \in \F$ such that 
$f \subset \Ba(\Mm(\F))(\epsilon)$\\
\hspace{1cm}if and only if\\
\>$\F \supset  \F \circ \Mm ( \F )$.
\end{tabbing}
\epf

\begin{fact}
\label{term} 
For any module $M: \xymatrix{I \ar[r]|{\circ} & A}$ if 
$- * M: [A,\R] \rightarrow \R$ preserves 
the terminal object (i.e. $\bigwedge_{ x \in A} M(x) = 0$)
then $\F(M)$ is a filter on $A$ with basis the family $\Ba(M)$.
\end{fact}
\pf  
Let us see first that the set of subsets of the form
$\Ba(M)(\epsilon)$ for $\epsilon > 0$, is a filter
basis on $A$.
Since $\bigwedge_{ x \in A} M(x) = 0$, for any $\epsilon > 0$
there is one $x$ with $M(x) \leq \epsilon$, i.e. 
$\Ba(M)(\epsilon) \neq \emptyset$. That
$\Ba(M)$ is a cofiltered subset of $\wp(A)$ ordered 
by inclusion is trivial. 
\epf

\begin{fact}
\label{counit}
If $M$ is a left module on $A$
then for all $x$, 
$$M(x) \leq \bigvee_{\epsilon > 0} \bigwedge_{y \mid M(y) 
\leq \epsilon} A(x,y).$$
\end{fact}
\pf
Let $x \in A$.
For all $y \in A$, $M(x) \leq M(y) + A(x,y)$ thus
for all $y \in A$, such that $M(y) \leq \alpha$, $M(x) \leq A(x,y) + \alpha$
and $M(x) \leq \bigwedge_{y \mid M(y) \leq \alpha} A(x,y) + \alpha$.
Consider $\epsilon > 0$. 
The map $\alpha \mapsto \bigwedge_{y \mid M(y) \leq \alpha} A(x,y)$
reverses the order so $\bigwedge_{ y \mid M(y) \leq \epsilon } A(x,y)$
$=$ $\bigvee_{ \alpha \geq \epsilon} \bigwedge_{ y \mid M(y) \leq \alpha } A(x,y)$
and
\begin{center}
 $(*)$
$M(x)$ $\leq$ $\bigvee_{ \alpha \geq \epsilon} \bigwedge_{ y \mid M(y)
  \leq \alpha } A(x,y) + \epsilon$.
\end{center}
Also for any $\alpha \leq \epsilon$,
\begin{tabbing}
\hspace{1cm}$M(x)$ \=$\leq$  \=$( \bigwedge_{ y \mid M(y) \leq \alpha } A(x,y) ) + \alpha$\\  
\>$\leq$  \>$( \bigwedge_{ y \mid M(y) \leq \alpha } A(x,y) ) + \epsilon$
\end{tabbing}
and thus 
\begin{center}
$(**)$
$M(x)$ $\leq$  $\bigvee_{\alpha \leq \epsilon}  \bigwedge_{ y \mid
  M(y) \leq \alpha } A(x,y)  + \epsilon$.
\end{center}
$(*)$ and $(**)$ give
$M(x) \leq  \bigvee_{\alpha > 0} \bigwedge_{y \mid M(y) \leq \alpha}A(x,y) + \epsilon$.
\epf

\begin{fact}
\label{P1flatM}
If the module $M: \xymatrix{I \ar[r]|{\circ} & A}$ 
is such that $-*M$ preserves cotensors
then 
$M(x) \geq \bigvee_{\alpha > 0} \bigwedge_{M(y) \leq \alpha} A(x,y)$.
\end{fact}
\pf
We show that for any $\epsilon > 0$ and any $x$ 
with $M(x) < \epsilon$ and any $\alpha > 0$,
there is a $y$ such $M(y) \leq \alpha$ and
$A(x,y) \leq \epsilon$.
Consider $\epsilon > 0$ and $x$ with
$M(x) < \epsilon$. Then 
\begin{tabbing}
\hspace{1cm}$0$ \=$=$ \=$[M(x), M(x)]$\\
\>$=$ \>$[M(x), A(x,-) * M ]$\\
\>$=$ \>$( M(x) \pitchfork A(x,-) ) * M$\\
\>$=$ \>$\bigwedge_{y \in A} ( M(y) + [M(x), A(x,y)] )$.
\end{tabbing}
So for any $\delta > 0$, there is a $y$ such that 
$M(y) + [M(x), A(x,y)] \leq \delta$.
This $y$ satisfies $M(y) \leq \delta$,
and $A(x,y) \leq M(x) + \delta$.
Now given any $\alpha>0$,
one may find a $y$ as required by
considering $\delta = min \{ \alpha, \epsilon - M(x) \}$.
\epf

According to \ref{counit} and \ref{P1flatM} above,
for any $\Pu$-flat module $M$, $M = \Mm \circ \F(M)$,
so that the filter $\F(M)$ is equal to 
$\F \circ \Mm \circ \F (M)$, hence is of type $1$
by \ref{unit}. This also shows that the functor
$\Mm: \FFu(A) \rightarrow \Pu^+(A)$ is essentially
surjective on objects, hence is part of an 
equivalence of category, with equivalence inverse
given by $M \mapsto \F(M)$.\\

To conclude the proof of Theorem \ref{mainres}, it remains to 
reformulate the distance on $\FFu(A)$ in purely metric terms. 
This is done in \ref{wfHom} below.
\begin{prop}
\label{dwflat2}
For any left module $M$ on $A$ and any filter $\F$,
$$[A^{op}, \R] ( \Mm(\F), M ) \leq \Limsup_{\F} M.$$
If $\F$ has type 1 then the inequality above
becomes an equality.
\end{prop}
\pf
To simplify notations, let $LHS$ and $RHS$ denote respectively  
$\bigvee_{x \in A} [\Mm(\F)(x),M(x)]$
and $\bigwedge_{f \in \F} \bigvee_{z \in f} M(z)$.\\

According to the definition of $\Mm(\F)$, for all $x \in A$,
for all $f \in \F$, $\Mm(\F)(x) \geq \bigwedge_{z \in f} A(x,z)$. 
So for any $x \in A$, any $f \in \F$ and any $\epsilon >0$, 
there exists a $z \in f$ such that 
$A(x,z) \leq \Mm(\F)(x) + \epsilon$.
For such a $z$, $M(x) \leq M(z) + A(x,z)$ and
$M(x) \leq M(z) + \Mm(\F)(x) + \epsilon$.
So, 
\begin{tabbing}
\hspace{1cm}\=$\forall x \in A, \forall f \in \F, \forall \epsilon >
0, \exists z \in f$,
$[\Mm(\F)(x),M(x)] \leq M(z) + \epsilon$,\\
thus \>
$\forall x \in A, \forall f \in \F, \forall \epsilon > 0$,
$[\Mm(\F)(x),M(x)] \leq (\bigvee_{z \in f} M(z)) + \epsilon$,\\
thus \>
$\forall f \in \F, \forall \epsilon > 0$,
$LHS \leq 
(\bigvee_{z \in f} M(z)) + \epsilon$,\\
thus \>
$\forall f \in \F$, 
$LHS \leq \bigvee_{z \in f} M(z)$,\\
thus \>
$LHS \leq RHS$.\\
\end{tabbing}
Suppose now that $\F$ has type 1.
Consider $\epsilon > 0$. One may find an $f_{\epsilon} \in \F$ 
such that for all $z \in f_{\epsilon}$, 
$\Mm(\F)(z) \leq \epsilon$. 
So,
\begin{tabbing}
\hspace{1cm}\=for any $z \in f_{\epsilon}$, $M(z) \leq [\Mm(\F)(z), M(z)] +
\epsilon$,\\
thus \>
$\bigvee_{z \in f_{\epsilon}} M(z) \leq 
(\bigvee_{x \in A}  [\Mm(\F)(x), M(x)]) + \epsilon$,\\
and \>
$RHS \leq LHS + \epsilon$.
\end{tabbing}
\epf
\begin{coro}
\label{wfHom}
For any filters $\F_1$ and $\F_2$ both of type 1,
$$[A^{op}, \R](\Mm(\F_1), \Mm(\F_2)) =  
\Limsup_{x \in \F_1} \Liminf_{y \in \F_2} A(x,y).$$ 
\end{coro}

A few remarks are in order.
First the underlying ordinary category $\FFu(A)_0$ 
of $\FFu(A)$, is a preorder, defined for all $\F_1$, $\F_2$ by 
\begin{fact}
$\F_1 \rightarrow \F_2$ $\Leftrightarrow$ $\F_1 \supset \F \circ \Mm(\F_2).$
\end{fact}  
\pf
\begin{tabbing}
\hspace{3cm}\=$\F_1 \rightarrow \F_2$\\
if and only if
\>$0 \geq [A^{op},\V](\Mm(\F_1),\Mm(\F_2))$\\
if and only if 
\>$\Limsup_{x \in \F_1} \Mm(\F_2)(x) = 0$\\
if and only if 
\>$\forall \epsilon > 0, \exists f \in \F_1, \forall x \in f, 
\Mm(\F_2)(x) \leq \epsilon$\\
if and only if 
\>$\forall \epsilon > 0, \exists f \in \F_1,
f \subset \Ba(\Mm(\F_2))(\epsilon)$\\
if and only if 
\>$\F_1 \supset \F \circ \Mm(\F_2)$.
\end{tabbing}
\epf
Note also that for any filters $\F_1$ and
$\F_2$, if $\F_1 \supset \F_2$ then 
$\Mm(\F_1) \Rightarrow \Mm(\F_2)$ and 
for any modules $M$, $N$, if $M \Rightarrow N$
then $\F(M) \supset \F(N)$. So that 
\ref{unit}, \ref{counit} and \ref{P1flatM} above shows 
the existence of an ordinary reflection of
the category of filters of type $1$
on $A$ with reverse inclusion ordering,
in the ordinary category of $\Pu$-flat modules on $A$. 
The ordinary category ${\FFu(A)}_0$ is exactly the 
category of fractions deduced from this reflection 
(see \cite{Bor94-1}, prop.5.3.1).

The following point results also from \ref{unit}.
The class of filters of type $1$ is the largest
class of filters such that fitted with the reverse 
inclusion ordering, it forms a category with a full 
reflection in the ordinary category of modules on $A$,
with functors given by the pair of maps
$M \mapsto \F(M)$ and $\F \mapsto \Mm(\F)$. 

Observe eventually that \ref{fac22} and \ref{dwflat2} 
express respectively the colimits and 
limits weighted by $\Pu$-flat modules $M$
in terms of the corresponding filters $\F(M)$.

\subsection{The theory of filters of type1}
Using the categorical machinery of section \ref{flatness}, 
we explore now the topology induced by the filters
of type 1. In particular we shall explicit in topological/metric 
terms the free $\Pu^+$-cocompletions.


One has a notion of non-symmetric convergence in $A$.
The {\em neighborhood filter} of $x \in A$, 
denoted $V_A(x)$, is 
the filter generated by the family of subsets  
$\{y \mid A(y,x) \leq \epsilon \}$ with $\epsilon > 0$.
Which is to say that $V_A(x)$ is $\F(A(-,x))$.
Given a filter $\F$ on $A$ and 
$x \in A$, we say that $\F$ {\em converges} to $x$,
that we write $\F \rightarrow x$, 
if and only if $\F \supset V_A(x)$.
If $\F$ has type 1 then   
$\F$ converges to $x$ if and only
if $\Mm(\F) \Rightarrow A(-,x)$.
By Yoneda, this is also equivalent to say 
that for any $a \in A$, 
$$A(x,a) \geq [A^{op},\R](\Mm(\F), A(-,a)),$$
or according to \ref{dwflat2} that 
$$A(x,a) \geq \Limsup_{\F} A(-,a).$$

\begin{defi}
A filter $\F$ on $A$ has {\em representative} $x_0$
if and only if for all $a \in A$,
$$A(x_0,a) = \Limsup_{\F}A(-,a).$$ 
\end{defi} 
Which is exactly to say that $x_0$ is the colimit $\Mm(\F)* 1$.
In particular if a representative of $\F$ exists
then it is unique up to isomorphism. In this case we denote it 
$rep(\F)$. Note that according to Yoneda, 
$rep(\F)$, when it exists, is necessarily the 
greatest lower bound in $A_0$ of the set of objects that $\F$ 
converges to.\\

Given a filter $\F$ on $A$ and a map $G:A \rightarrow B$
the direct image of $\F$ denoted $G(\F)$ is the filter
on $B$ generated by the family of subsets $G(f)$ for $f \in \F$.
It is easy to check for $\F$ and $G$ as above that 
if $G$ is non-expansive
and $\F$ has type $1$
then $G(\F)$ has again type $1$. Moreover,
\begin{prop}
\label{image}
Given a filter $\F$ of type 1 on $A$, and a functor $G: A \rightarrow B$,
$\Mm(G(\F)): B^{op} \rightarrow \R$ is the (pointwise) left Kan extension of 
$\Mm(\F): A^{op} \rightarrow \R$ along $G^{op}$.
\end{prop}
\pf One has the pointwise computation (see \cite{Kel82},
(4.17), p.115), 
\begin{tabbing}
$Lan_{G^{op}}(\Mm(\F))(b)$ \=$=$ \=$B(b,G-) * \Mm(\F)$\\
\>$=$ \>$\bigvee_{f \in \F} \bigwedge_{x \in f} B(b,Gx)$, 
according to \ref{fac22},\\
\>$=$ \>$\bigvee_{g \in G(\F)} \bigwedge_{y \in g} B(b,y)$\\
\>$=$ \>$\Mm(G(\F))(b)$.
\end{tabbing}
\epf
Note that we could already infer from \ref{Lprpf}
that for any filter $\F$ and any functor $G: A \rightarrow B$,
if $\Mm(\F)$ is $\Pu$-flat then  $Lan_{G^{op}}(\Mm(\F))$
is also $\Pu$-flat.
Now remark that for any weight $F:A^{op} \rightarrow \V$,
for any $G:A \rightarrow B$ and any $H: B \rightarrow C$ 
with $B$ and $C$ small, $Lan_{G^{op}}(F)*H \cong F*HG$.
So that one has the following consequence of \ref{image}.
\begin{fact} 
\label{trans}
Given a filter $\F$ of type 1 on $A$ and
a non expansive map $G: A \rightarrow B$,
$\Mm(\F) * G$ is 
(up to isomorphism) the representative of
the filter $G(\F)$ of type $1$.
\end{fact}

We shall call a general metric space $A$ {\em (type 1)-complete}, 
when any filter of type 1
on $A$ admits a representative, according to \ref{trans},
this is to say when the associated category is
{\em $\Pu^+$-cocomplete}. 
Now consider a filter $\F$ of type 1 on $K$, 
and non-expansive maps $G:K \rightarrow A$, 
and a functor $H:A \rightarrow B$. 
Then $H$ (as a functor) preserves the 
colimit $\Mm(\F)*G$ if and only $H$ (as a non-expansive map) 
preserves the representative of $G(\F)$, i.e.
$$H( rep(G(\F)) ) = rep (H \circ G( \F  ) ).$$
To sum up: the $\R$-functors preserving the $\Pu^+$-colimits 
are exactly 
the non-expansive maps preserving the representatives of
filters of type $1$.
Now a translation of Kelly's theorem (\ref{freecoc}) gives
the following completion for metric spaces.
\begin{theo}\label{type1comp}
For any general metric space $A$, there exists 
a (type $1$)-complete metric space 
$\bar{A}$ together with a non expansive map 
$i_A: A \rightarrow \bar{A}$, such that 
for any (type $1$)-complete general metric space $B$,
composing with $i_A$ defines an equivalence of general 
metric space:
\begin{center}
$- \circ i_A : [\bar{A},B]' \cong [A,B]$
\end{center}
where
$[\bar{A},B]'$ denotes the sub-metric space of 
$[\bar{A},B]$ of maps preserving the representatives
of filters of type $1$.\\
Any $\bar{A}$ as above is equivalent to the metric space $\FFu(A)$,
and the embedding $i_A$ is, up to this 
equivalence, the non-expansive map $A \rightarrow \FFu(A)$ sending
a point $x$ to its neighborhood filter $V_A(x)$.
\end{theo} 
Let us call therefore $\FFu(A)$ the 
{\em completion of type 1} of $A$.
Let us add that with $A$ and $B$ as above, 
for any $f: A \rightarrow B$, 
the unique extension $\bar{f}: \FFu(A) \rightarrow B$
of $f$ through $i_A$ sends any filter $\F$ to
the representative of its direct image by $f$ in $B$.
To check this, just come back to the categorical formulation.
From \cite{Kel82} Theorem 4.97,  
$\bar{f}$ is the left Kan extension of $f$ along $i_A$
and sends any $M$ in $\Pu^+(A)$ to $M*f$.
Translate then using \ref{trans}.\\

Let us mention the following result that can be inferred from
purely categorical arguments. The limits in \ref{wfHom} 
``commute'' when the first argument is Cauchy:
\begin{fact}
\label{CcHom}
For any Cauchy filter $\F_1$ and any filter $\F_2$ of type 1,
$$[A^{op}, \R](\Mm(\F_1),\Mm( \F_2)) = \Liminf_{y \in \F_2} \Limsup_{x \in \F_1} A(x,y).$$
\end{fact}
To see this we shall need the following result
see for instance \cite{Bet85}-remark 4 p.171
or \cite{KS05} section 6.
\begin{prop}
\label{Homladj}
Given a complete and cocomplete monoidal closed $\V$, 
if a $\V$-module $M: \xymatrix{I \ar[r]|{\circ} & C}$
has a right adjoint $\tilde{M}: \xymatrix{C \ar[r]|{\circ} & I}$
then for any $\V$-module $N:\xymatrix{I \ar[r]|{\circ} & C}$, 
$[C^{op}, \V](M,N) \cong \tilde{M} * N$.
\end{prop}

We can now prove \ref{CcHom}.\\ 
\pf 
Recall that for any Cauchy filter $\F$ on $A$,
$\Mm(\F)$ $=$ $\Mp(\F)$ $=$ $\Ml(\F)$ 
and this left module on $A$ has right adjoint the module 
$\Mr(\F)$ defined 
by the map $x \mapsto \Limsup_{y \in \F} A(y,x) = \Liminf_{y \in f} A(y,x)$.
So according to \ref{Homladj} and \ref{fac22}, for any Cauchy filter $\F_1$ and any 
filter $\F_2$ of type 1,
\begin{tabbing}
\hspace{1cm}$[A^{op}, \R]( \Mm(\F_1), \Mm(\F_2) )$ \=$=$  \=$\Mr(\F_1) * \Mm(\F_2)$\\
\>$=$
\>$\Liminf_{y \in \F_2} \Mr(\F_1)(y)$\\ 
\>$=$
\>$\Liminf_{y \in \F_2} \Limsup_{x \in \F_1} A(x,y).$
\end{tabbing}
\epf


\subsection{$\Pa$-flat modules and filters of type $\aleph$.}
$\aleph$ will denote in this section any infinite regular cardinal.
We turn now to the case of $\Pa$-flat modules
and their corresponding filters, the so-called
filters of type $\aleph$.

\begin{defi}
A filter $\F$ on $A$ has {\em type $\aleph$}
if and only if for any $\epsilon > 0$, there exists an $f \in \F$
such that for any family of elements ${(x_i)}_{i \in I}$ 
of $f$, with $\sharp I < \aleph$, for any $g \in \F$, 
there exists $y \in g$ such that $A(x_i,y) \leq \epsilon$.\\
\end{defi}

Note the inclusion of classes of filters:
\begin{fact}
Cauchy $\Rightarrow$ type $\aleph$ $\Rightarrow$ 
type $\omega$ $\Rightarrow$ type $1$.
\end{fact}
Also when $A$ is symmetric, that is when 
$A(x,y) = A(y,x)$, filters of type $\omega$ are
also Cauchy. (We shall see later a few consequences 
of this fact.)\\

With the above definition of filters, the equivalence of Theorem 
\ref{mainres} restricts to the full subcategories
of $\Pa$-flat presheaves and filters of 
type $\aleph$.
\begin{theo}\label{folk}
The full subcategory of $\FFu(A)$ induced
by the filters of type $\aleph$ is equivalent
to $\Pa^+(A)$, the equivalence functors being
given by the maps $\F \mapsto \Mm(\F)$ and
$M \mapsto \F(M)$.  
\end{theo}
This results from Theorem \ref{mainres} and
\ref{flatf} and \ref{lex} below.

\begin{fact}
\label{flatf} 
For any $\Pa$-flat module $M: \xymatrix{\I \ar[r]|{\circ} & A}$,   
$\F(M)$ has type $\aleph$.
\end{fact}
\pf
If $- * M$ preserves conical $\aleph$-limits then it preserves
in particular the terminal object 
and according to \ref{term}, $\F(M)$ is a filter on $A$. 
The fact that the filter basis $\Ba(M)$ 
generates a filter of type $\aleph$ is a consequence of the following 
result.
\epf

\begin{fact}
\label{P2flatM}
If $M: \xymatrix{\I \ar[r]|{\circ} & A}$ is $\Pa$-flat 
then for any $\epsilon > 0$ and
any family ${(x_i)}_{i \in I}$, with $\sharp I < \aleph$,
such that for all $i$, $M(x_i) \leq \epsilon/2$ and any $\alpha > 0$,
there is a $y$ such that $M(y) \leq \alpha$ and for all $i \in I$,
$A(x_i,y) \leq \epsilon$.
\end{fact}
\pf
$-*M$ preserves conical $\aleph$-limits and cotensors.
Consider $\epsilon > 0$ and a family $(x_i)_{i \in I}$'s 
with $\sharp I < \aleph$  such that 
$M(x_i) \leq \epsilon / 2$. 
Let us write $\epsilon' = \bigvee_{i \in I}M(x_i)$. Then 
\begin{tabbing}
\hspace{1cm}$0$ \=$=$ \=$[\epsilon', \bigvee_{i \in I}M(x_i)]$\\
\>$=$ 
\>$[\epsilon', \bigvee_{i \in I}( A(x_i,-) * M )]$\\
\>$=$
\>$[\epsilon', ( \bigvee_{i \in I}A(x_i,-) ) * M ]$\\ 
\>$=$$(  \epsilon' \pitchfork (\bigvee_{i \in I}A(x_i,-) ) * M$\\
\>$=$$\bigwedge_{y \in A} (M(y) + [\epsilon', \bigvee_{i \in I}A(x_i,y)])$.
\end{tabbing}
So for any $\delta > 0$, there is a $y$ such that 
$M(y) + [\epsilon', \bigvee_{i \in I}A(x_i,y)] \leq \delta$.
This $y$ satisfies $M(y) \leq \delta$,
and for all $i$, $A(x_i,y) \leq \epsilon' + \delta$.
Now given any $\alpha > 0$,
one may find a $y$ as required
 by considering 
$\delta = min \{ \alpha, \epsilon - \epsilon' \}$.
\epf

\begin{fact}
\label{lex}
If the filter $\F$ on $A$ has type $\aleph$ then
$-* \Mm(\F)$ preserves conical $\aleph$-limits, i.e.
for any family $(N_i)_{i \in I}$
of right modules on $A$, with $\sharp I < \aleph$,
$$\bigwedge_{x \in A} ( \Mm(\F)(x) + \bigvee_{i \in I} N_i(x)  ) =  
\bigvee_{i \in I} \bigwedge_{x \in A} ( \Mm(\F)(x) +  N_i(x) ).$$
\end{fact}
\pf We only need to prove 
$$\bigwedge_{x \in A} ( \Mm(\F)(x) + \bigvee_{i \in I} N_i(x)  ) 
\leq  \bigvee_{i \in I} \bigwedge_{x \in A} ( \Mm(\F)(x) +  N_i(x) ).$$

Let $\epsilon > 0$.
If there is a filter $\F$ on $A$ then $A$ is not empty and
for each $i \in I$, there is an $x_i \in A$ such that
$$N_i * \Mm(\F) + \epsilon 
= \bigwedge_{x \in A} (\Mm(\F)(x) + N_i(x)) + \epsilon \geq 
\Mm(\F)(x_i) + N_i(x_i).$$
Let $f \in \F$. Given a family
of $x_i$'s as above, for each $i$,
$\Mm(\F)(x_i) \geq \bigwedge_{y \in f} A(x_i,y)$, 
thus there is an $y_i \in f$ such that
$\Mm(\F)(x_i) + \epsilon \geq A(x_i, y_i)$ and
\begin{tabbing}
\hspace{1cm}
$2 \cdot \epsilon + N_i * \Mm(\F)$  \=$\geq$ 
\=$A(x_i, y_i) + N_i(x_i)$\\
\>$\geq$ \>$N_i(y_i)$.
\end{tabbing}
Since $\F$ has type $\aleph$, we can choose $f$ so that  
for the $y_i \in f$ as above, for all $g \in \F$, there exists 
$z \in g$ such that for all $i$, $A(y_i,z) \leq \epsilon$.\\

Thus for all $\epsilon >0$,
for all $g \in \F$, there exists $z \in g$ such that 
for all $i \in I$, 
\begin{tabbing}
\hspace{1cm}
$3 \cdot \epsilon   + N_i * \Mm(\F)$
\=$\geq$   
\=$A(y_i,z) + N_i(y_i)$ for some suitable $y_i$'s,\\
\>$\geq$ \>$N_i(z)$.
\end{tabbing}
and thus 
\begin{tabbing}
\hspace{1cm}\=$\forall \epsilon >0, \forall g \in \F, \exists z \in g$, 
$\epsilon + \bigvee_{i \in I}( N_i * \Mm(\F) ) \geq 
\bigvee_{i \in I} N_i(z);$\\
thus \>$\forall \epsilon >0, \forall g \in \F$, 
$\epsilon + \bigvee_{i \in I}( N_i * \Mm(\F) ) \geq 
\bigwedge_{z \in g}\bigvee_{i \in I} N_i(z);$\\
thus \>$\forall g \in \F, 
\bigvee_{i \in I}( N_i * \Mm(\F) ) \geq 
\bigwedge_{z \in g} \bigvee_{i \in I} N_i(z).$
\end{tabbing}
So \begin{tabbing}
\hspace{1cm}$\bigvee_{i \in I} ( N_i * \Mm(\F) )$ \=$\geq$
\=$\bigvee_{g \in \F} \bigwedge_{z \in g} \bigvee_{i \in I} N_i(z)$\\
\>$=$ \>$( \bigvee_{i \in I} N_i ) * \Mm(\F)$, according to \ref{fac22}. 
\end{tabbing}
\epf

One obtains also (type $\aleph$)-completions for general metric 
spaces. Let $\FFa(A)$ stand for the full subcategories
of $\FFu(A)$ generated by the filters of type $\aleph$.
Then Theorem \ref{type1comp} still 
holds after that ``type $\aleph$'' and ``$\FFa(A)$''
has been substituted everywhere respectively 
to ``type $1$'' and ``$\FFu(A)$''.
Actually the only point to check to establish this,
is that the direct image by a non expansive map
of a filter a type $\aleph$, is again of type $\aleph$;
which is straightforward.

\subsection{Examples}
We give now examples of filters of types $1$ and 
$\omega$, complete spaces and completions.\\ 


Recall from \cite{Kel82} (3.74) that
any monoidal closed $\V$ that is complete and cocomplete
as an ordinary category, is complete and cocomplete as a 
$\V$-category. Thus
\begin{fact}
$\R$ is (type $1$)-complete.
\end{fact}
 
Since $\Pu \supset \C$ one has $\Pu^+ \subset \C^+$ and thus
according to \ref{flacot}, one has also the following.
\begin{prop}
The (type $1$)-completion of $\R$ is equivalent to $\R$.
\end{prop}
It might be worth detailing a bit the situation here.
For a filter $\F$ on $\R$, we write $\liminf(\F)$ 
for $\Liminf_{\F} id$ which is $\bigvee_{f \in \F} \bigwedge_{x \in f} x$.
If $\F$ has type $1$ then $\liminf(\F)$ is just $\Mm(\F)*1$
by Proposition \ref{fac22},
and according to Proposition \ref{flacot}, 
$\Mm(\F): \V^{op} \rightarrow \V$ is isomorphic to $[-, \liminf(\F)]$. 
In the latter case,
$\F \supset V_{\cal V}(\liminf \F)$ i.e. $\F$ converges
to $\liminf(\F)$.

\begin{prop}
\label{wfeqf}
Filters of type $1$ on $\R$ have type $\omega$.
\end{prop}
\pf 
Let $\F$ be a filter of type $1$ on $\R$.
Consider $\epsilon > 0$ then there exists an $f \in \F$ 
such that for any $x \in f$, for any $g \in \F$, there 
exists a $y$ such that $A(x,y) \leq \epsilon$.
For this $f$, for any finite family of elements $x_i$ in $\F$ 
and  for any $g \in \F$, one may find elements $y_i \in g$,
such that for all $i$, 
$[x_i, y_i] \leq \epsilon$, i.e. 
$y_i \leq x_i + \epsilon$. Choosing the least of those $y_i$'s, say 
$z$, one has $[x_i,z] \leq \epsilon$ for all $i$. 
\epf

We investigate now the case of {\em symmetric} spaces. 
\begin{fact}
\label{symA}
If $A$ is symmetric, 
\begin{itemize}  
\item $(1)$ filters of type $\omega$ on $A$ are Cauchy;
\item $(2)$ For any Cauchy filter $\F$, $\Ml(\F)$ and $\Mr(\F)$ 
have the same underlying map. 
\item $(3)$ Any left adjoint module on $A$ 
has the same underlying map as its right adjoint;
\item $(4)$ For any left adjoint module $M$ on $A$,
$\F(M) = \Fs(M)$;
\item $(5)$ $\Pd$-flat modules are left adjoint;
\item $(6)$ The ordinary category ${\Q(A)}_0$ is discrete.
\end{itemize}
\end{fact}
\pf
$(1)$ -- that was already mentioned -- and $(2)$ are trivial.
$(3)$ holds since for any left module $M$ with right adjoint $\tilde{M}$,
according to $(2)$ their 
underlying maps satisfy $M = \Ml \circ \Fs(M) = \Mr \circ \Fs(M) = \raM$.
$(4)$ is straightforward from $(3)$. To prove  
$(5)$, consider the successive equivalences
\begin{tabbing}
\hspace{1cm}\=$M$ is $\Pd$-flat\\
\>if and only if
$M = \Mm(\F)$ for a filter of type $\omega$\\
\>if and only if
$M = \Mm(\F)$ for a Cauchy filter (according to $(1)$),\\
\>if and only if
$M$ is left adjoint.
\end{tabbing}
Now we show $(6)$, namely: the underlying
subcategory of the full subcategory of presheaves
$[A^{op},\R]$ with objects left adjoint modules is
discrete
(in the particular case $\V = \R$). 
For any left adjoint module $M$ on $A$, 
$M$ has the same underlying map as
its right adjoint $\raM$.
Now consider another left adjoint module $N$ on $A$,
with right adjoint $\tilde{N}$.
Then $M \Rightarrow N$ if and only if
$\forall x \in A, M(x) \geq N(x)$ 
if and only if
$\forall x \in A, \tilde{M}(x) \geq \tilde{N}(x)$
if and only if
$\tilde{M} \Rightarrow \tilde{N}$.
But also if $M \Rightarrow N$ then 
$1 \Rightarrow \tilde{M}N$
since $M \dashv \tilde{M}$
and then
$\tilde{N} \Rightarrow \tilde{M}$
since $N \dashv \tilde{N}$.  
So $M \Rightarrow N$ if and only if $M=N$.
\epf

The assertion $(1)$ above tells us that 
when the general metric space $A$ is symmetric,
it happens that
$\Pd^+[A] = \p^+[A]$, which is also to say
that the completion of type $\omega$ of $A$ 
is its Cauchy completion. 
Nevertheless even when $A$ is symmetric, its completion 
of type $1$ may be not symmetric. We show below that 
it consists of the set of non-empty closed subsets of $A$ 
with what one could call a ``semi-Hausdorff'' distance. 

\begin{prop}
\label{sycomp}
The completion of type $1$ of a symmetric $A$ is the set 
of non-empty closed subsets of its Cauchy-completion $\bar{A}$
with pseudo distance $d$ given by 
$d(X,Y) = \bigvee_{x \in X} \bigwedge_{y \in Y} \bar{A}(x,y)$. 
\end{prop}

To prove this result, we shall establish 
a characterization of filters
of type $1$ as certain colimits of {\em forward Cauchy sequences}.
These sequences belong to folklore and 
were introduced as a generalization of the classical 
Cauchy sequences for non-symmetric metric spaces. It 
is not known by the author whether they admit a reasonable description 
in categorical terms. Nevertheless these sequences define peculiar 
filters of type $\omega$ and as such, their whole class has a surprising 
density property in ${\FFu(A)}_0$ (see
Theorem \ref{charffil} below).\\
 
Given a sequence $(x_n)_{n \in \nit}$ on $A$, the associated filter,
still denoted $(x_n)$,   
has basis the family of sets $\{ x_p \mid p \geq n \}$.
We say that $(x_n)$ is:\\
- {\em of type 1}, respectively {\em of type $\omega$}, if the 
associated filter is so;\\
- {\em forward Cauchy} if and only
if $\forall \epsilon > 0, \exists N \in \nit, 
\forall m \geq n \geq N, A(x_n,x_m) \leq \epsilon$.\\
Note that any forward Cauchy sequence is obviously of type $\omega$.

\begin{theo}
\label{charffil}
${\FFu(A)}_0$ has all the colimits of non-empty diagrams
and filters of type $1$ are 
colimits in ${\FFu(A)}_0$ of non-empty diagrams
with values forward Cauchy sequences. 
\end{theo}
This result will follow from \ref{yogl3} and \ref{yogl1} below.
We note first that the base $\R$ has a very peculiar property:
\begin{fact}
The conical terminal object commutes with 
conical colimits of non-empty diagrams in $\R$.
\end{fact}
\pf
The weight for the conical terminal object is
the unique functor $!: \emptyset \rightarrow \R$, where
$\emptyset$ denotes the empty category. 
Also the functor category $[\emptyset,\R]$ is 
isomorphic to the terminal category $1$ with one object $*$ and 
hom $1(*,*) = 0$. Therefore the fact that for any $v \in \R$,
$[v,0] = 0$ shows that the limit $\{!,! \}$ is $0$. 
Now we consider an ordinary category $J$ and let $J^{\sharp}$
denote the free $\R$-category over $J$.
Since limits in functor categories are pointwise,
the limit weighted by $!$ of the unique 
functor $\emptyset \rightarrow [J^{\sharp},\R]$,
denoted again $\{!,!\}$, is the constant functor $J^{\sharp}
\rightarrow \R$ to $0$. The conical colimit of this functor is 
$0$ since $J$ is non-empty.
On the other hand,
the unique functor $J^{\sharp} \rightarrow [\emptyset,\R] 
\cong 1$ is necessarily the
constant one with image the unique $!: \emptyset \rightarrow \R$,
its conical colimit is necessarily the functor 
$!:\emptyset \rightarrow \R$ which
limit weighted by $!: \emptyset \rightarrow \R$ is $\{!,!\}$,
which is known to be $0$. 
\epf

According to this and the fact that conical colimits 
of non-empty diagrams commute
also in $\R$ with cotensors (\ref{facR}), one has:
\begin{fact} 
the class of weights for conical colimits of non-empty
diagrams 
is contained in $\Pu^+$. 
\end{fact} 
Therefore, since $\Pu^+$ is saturated, one has the following result:
\begin{fact}
\label{yogl3} $\Pu^+(A)$ is closed in $[A^{op},\R]$ 
under the formation of conical colimits of non empty 
diagrams.
\end{fact}

\begin{fact}
\label{yogl1}
Given a filter $\F$ of type $1$ on $A$ and a left module $M$
such that $\Mm(\F) \not \Rightarrow M$
(i.e. $\Mm(\F) \not \geq M$), 
there is a forward Cauchy sequence 
$(y_n)$ such that $(y_n) \rightarrow \F$ and
$\Mm(y_n) \not \Rightarrow M$.
\end{fact}
\pf
The arguments will be technical and we introduce some convenient  
notation. For any $f \subset A$, any $\epsilon > 0$ and, 
any $F \subset \wp\wp(A)$, 
we let $P(f,\epsilon, F)$ denote the property:\\
``for all $x$ in $f$, for 
all $g \in F$ there
exists $y \in g$ such that
$A(x,y) \leq \epsilon$''.\\
 
Now, by hypothesis 
there exists $x \in A$ such that 
$\bigvee_{f \in \F} \bigwedge_{y \in f} A(x,y) < M(x)$. 
Consider such an $x$. There exists $\alpha> 0$ 
such that for any $f \in \F$, there exists $y \in f$ such that 
$A(x,y) + \alpha < M(x)$. Note then that for such a $y$,
$A(x,y)$ is necessarily finite.\\

Since $\F$ has type 1, one can define
a sequence $(f_n)$ of elements of $\F$ such that for all $n \in \nit$, 
$f_{n+1} \subset f_n$ and $P(f_n , \alpha \cdot 2^{-2-n}, \F)$.
$(f_n)$ is defined inductively as follows.\\

Choose first $f_0$ such that $P(f_0, \alpha \cdot 2^{-2}, \F)$.\\

If $f_n$ is defined then 
one can find $g \in \F$ such that $P(g, \alpha \cdot 2^{-2-(n+1)}, \F)$
and let $f_{n+1} = f_n \cap g$.\\

Then one can build a sequence $(y_n)$ where
for all integer $n$, $y_n \in f_n$,
$y_0$ is such that $A(x,y_0) + \alpha < M(x)$,
and for all integer $n$, $y_n \in f_n$, $A(y_n,y_{n+1}) \leq \alpha \cdot 2^{-2-n}$.\\

Actually this ensures that:\\
$(1)$ $(y_n)$ is forward Cauchy;\\
$(2)$ $(y_n) \rightarrow \F$;\\
$(3)$ $\Mm(y_n) \not \Rightarrow M$.\\
  
$(1)$ holds since
for all $n \leq p \in \nit$,
\begin{tabbing}
\hspace{1cm}$A(y_n, y_p)$ \=$\leq$ \=$A(y_n, y_{n+1}) + ... + A(y_{p-1},y_p)$\\ 
\>$\leq$ \>$\alpha \cdot ( 2^{-2-n} + 2^{-2-(n+1)} + ... )$\\
\>$=$ \>$\alpha \cdot 2^{-1-n}$.
\end{tabbing}

$(2)$ holds since $(y_n)$ is forward Cauchy, 
for any $n \in \nit$, $\{ y_p / p \geq n \} \subset f_n$ and
$P(f_n, \alpha \cdot 2^{-2 -n},\F)$.\\ 

$(3)$ holds since
for all $n \in \nit$, 
\begin{tabbing}
\hspace{1cm}$A(x,y_n)$\=$\leq$\=$A(x,y_0) + A(y_0, y_1) + ... + A(y_{n-1}, y_{n})$\\
\>$\leq$\= $A(x,y_o) + \alpha/2$
\end{tabbing}
so $A(x,y_n) + \alpha/2 <  M(x)$.
Thus $\Mm(y_n)(x) = \bigvee_{n \in \nit} \bigwedge_{p \geq n} A(x,y_p) < M(x)$.
\epf

To finish proving Theorem \ref{charffil}, it remains to see 
that any filter $\F$ of type $1$ dominates
at least one forward Cauchy sequence. But this holds for such an
$\F$ according 
to \ref{yogl1}, since $\Limsup_{\F} \Mm(\F) = 0$ and 
thus $\Mm(\F) \not \Rightarrow +\infty$ where $+\infty$ denotes
here the constant module with value $+\infty$.\\

We come back now to the proof of Proposition \ref{sycomp}.\\ 
\pf
Since $A$ is symmetric, its Cauchy completion is equivalent to the 
metric space, say $B$, with objects Cauchy filters on $A$ with
``symmetric'' distance $d$ given for all $\varphi$, $\psi$ by
$d(\varphi, \psi)$ $=$ 
$[A^{op},\R](\Mm(\varphi), \Mm(\psi))$ $=$ 
$\Mr(\varphi) * \Mm(\psi)$
by Proposition \ref{Homladj}.
This metric space is also symmetric (according to \ref{symA})
and forward Cauchy sequences in $B$ are just the Cauchy ones.

Consider a filter $\F$ of type $1$ on $A$ and let $\bar{\F}$ 
denote the set of Cauchy filters $\phi$ such that 
$\phi \rightarrow \F$. According to
Theorem \ref{charffil} and \ref{yogl3}, $\bar{\F}$ is 
not empty and $\Mm(\F)$ is the pointwise conical colimit 
in $[A^{op},\V]$: 
$$\Mm(\F) = \bigwedge_{\varphi \in \bar{\F}} \Mm(\varphi).$$
 
Now for any subset $X$ of $B$ that satisfies the 
property $$(*)\;\;\Mm(\F) = \bigwedge_{\varphi \in X} \Mm(\varphi),$$
and any Cauchy filter $\psi$ on $A$, one has:
\begin{tabbing}
$d(\psi,\F)$ \=$=$ \=$[A^{op},\R](\Mm(\psi), \Mm(\F))$\\
                 \>$=$ \>$[A^{op},\R](\Mm(\psi),  \bigwedge_{\varphi
                 \in X}  \Mm(\varphi))$\\
\>$=$ \>$\bigwedge_{\varphi \in X} [A^{op},\R](\Mm(\psi),\Mm(\varphi))$
$(**)$\\
\>$=$ \>$\bigwedge_{\varphi \in X} d(\psi,\varphi).$
\end{tabbing}
where $(**)$ above holds since the presheaf $\Mm(\psi): 
A^{op} \rightarrow \R$ is a small projective (or equivalently 
the module $\Mm(\psi): \xymatrix{I \ar[r]|{\circ} & A}$ is 
left adjoint).


As a consequence of this, one has that for any subset $X$ of $B$
satisfying $(*)$, the adherence $\bar{X}$ of $X$ in $B$ is $\bar{\F}$.
Therefore $\bar{\F}$ is the only closed subset $X$ in $B$ 
satisfying the condition $(*)$. 
This result, together with Theorems \ref{mainres} and
\ref{charffil}, show that 
the equation $(*)$ defines a bijection between
$\Pu$-flat modules on $A$ and 
and non-empty subsets of $B$.

Eventually given two filters on $A$ of type $1$,
$\F_1$ and $\F_2$,
\begin{tabbing}
$[A^{op}, \R](\Mm(\F_1),\Mm(\F_2))$
\=$=$ \=$[A^{op}, \R](\bigwedge_{\varphi \in \bar{\F_1}}
\Mm(\varphi),\Mm(\F_2))$\\
\>$=$ \>$\bigvee_{\varphi \in \bar{\F_1}}[A^{op}, \R](\Mm(\varphi),
\Mm(\F_2))$\\
\>$=$ \>$\bigvee_{\varphi \in \bar{\F_1}}[A^{op}, \R](\Mm(\varphi),
\bigwedge_{\psi \in \bar{\F_2}} \Mm(\psi)  )$\\
\>$=$ \>$\bigvee_{ \varphi \in \bar{\F_1} } \bigwedge_{\psi \in \bar{\F_2}}
d(\varphi,\psi)$,
\end{tabbing}
since, again, $\Mm(\varphi)$
is a small projective
for any Cauchy filter $\varphi$.
\epf

\section{The case $\V = \TWO$.}\label{preorder}
Preorders as enrichments over the category $\TWO$
were mentioned in \cite{Law73}. After a brief reminder,
we shall characterize in the context $\V = \TWO$ 
the $\Pu$- and $\Pa$- flatnesses and the associated
free cocompletions.
We shall see that the free $\Pd^+$-cocompletion 
is the classic {\em dcpo} completion.\\

$\TWO$ stands for the two-object category generated by the graph
$\xymatrix{ 0 \ar[r] & 1}$. It is a partial order and has 
a monoidal structure with tensor $\wedge$
(the logical ``and'') and unit 
$1$. $\TWO$ is closed since for all $x,y,z \in \TWO$,
$$x \wedge y \leq z \Leftrightarrow x \leq (y \Rightarrow z)$$  
where $\Rightarrow$ denotes the usual entailment relation.
Small $\TWO$-categories are just preorders: for any 
small $\TWO$-category $A$, 
its associated preorder is defined by $x \rightarrow y$ if and only 
if $A(x,y) = 1$. Along the same line there is a bijection
between $\TWO$-functors and monotonous maps.
Any $\TWO$-module $M:\xymatrix{\I \ar[r]|{\circ} & A}$ 
corresponds to a downset $\cI_M = \{x \mid M(x) = 1 \}$
on the preorder $A$, and this correspondence between
modules and downsets is bijective.
Also for any preorders $A$ and $B$, the set of 
monotonous maps is considered as pointwise ordered, 
which corresponds to the $\TWO$-enriched 
categorical structure of the functor category $[A,B]$.
Recall that small $\TWO$-categories are always 
Cauchy-complete.
Also the {\em downward completion} of a preorder
may be described as its set of downsets ordered by inclusion, 
and is its free completion.
Via the above translations, the downward completion
is just the free cocompletion of $\TWO$-categories.\\

Let us turn now to the $\Pu$- and $\Pa$-flatnesses.
In the rest of this section, $A$ denotes a small 
$\TWO$-category that 
we freely consider as a preorder, and $\aleph$ is any
infinite regular cardinal. 
Using \ref{cotpoint}, \ref{precon} and \ref{colimpre} again,
one gets the following.
\begin{fact}
For any module $M: \xymatrix{I \ar[r]|{\circ} & A}$,
\begin{itemize}
\item $-*M: [A,\TWO] \rightarrow \TWO$ preserves
the (conical) terminal object i.e.
$1 * M = 1$, if and only if 
$$(1)\;\;\bigvee_{x \in A} M(x) = 1;$$  
\item $-*M$ preserves conical $\aleph$-limits if and only if\\ 
$(2)$ For any family of right modules 
$N_i: \xymatrix{A \ar[r]|{\circ} & I}$, $i$ ranging in $I$ and
$\sharp I < \aleph$,
$$\bigvee_{x \in A}( M(x) \wedge \bigwedge_{i \in I} N_i(x) )
= \bigwedge_{i \in I}( \bigvee_{x \in A} M(x) \wedge N_i(x)  );$$
\item $-*M$ preserves cotensors if and only if\\
$(3)$ For any $v \in \TWO$ and 
any right module $N: \xymatrix{A \ar[r]|{\circ} & I}$, 
$$\bigvee_{x \in A}( M(x) \wedge  (v \Rightarrow N(x)) ) 
\;\;=\;\; (\; v \Rightarrow  \bigvee_{x \in A}( M(x) \wedge N(x) )\; ).$$
\end{itemize}
\end{fact}
Condition $(1)$ above is equivalent to the fact that 
$\cI_M$
is not empty.
Condition $(3)$ reduces for $v=1$ to the trivial equation $N*M = N*M$. 
For $v= 0$, it reduces to $1 = \bigvee_{x \in A} M(x)$,
that is $(1)$ again.
Recall that a downset $\cI$ of $A$ is said {\em $\aleph$-directed}
if and only if it satisfies the property:\\
$(*)$ Any subset of $\cI$ of cardinality strictly less than $\aleph$
has an upper bound in $\cI$.\\
Note that such an $\I$ is necessarily not empty.
Now observe that condition $(2)$ above is equivalent to the fact that 
$\cI_M$ satisfies $(*)$. 
So one has a bijection between the following on $A$:
\begin{itemize}
\item $\Pu$-flat left modules and non-empty downsets, 
\item $\Pa$-flat left modules and $\aleph$-directed 
downsets.
\end{itemize}
Thus $\Pu^+(A)$ is equivalent as a preorder to
the set of non-empty downsets of $A$ with inclusion 
ordering, whereas $\Pa^+(A)$ is equivalent to the set the 
$\aleph$-directed downsets of $A$
with inclusion ordering, that we shall write $\adcpo(A)$.\\

Eventually given a weight $M: A^{op} \rightarrow \TWO$ and 
$G: A \rightarrow B$ with $B$ small, 
$b \in B$ is the colimit $M*G$ if and only 
if $b$ is the least upper bound in the preorder $B$ of the downset 
generated by the direct image of $\cI_M$ by $G$.
Also a functor $H: B \rightarrow C$ preserves
$M*G$ as above if and only if the corresponding
monotonous map preserves the least upper bound 
of $G(\cI_M)$.
Remember that a partial order is said {\em $\aleph$-directed complete}
or, for short, is a {\em $\aleph$-dcpo}, when it admits all least upper 
bounds for $\aleph$-directed subsets.
$\aleph$-dcpos with maps 
preserving $\aleph$-directed least upper bounds 
form the category $\adcpo$.
So we know from categorical considerations that $\adcpo(A)$
is $\aleph$-directed complete.\\

From \ref{freecoc} again, one can deduce a few completions for 
preorders, and in particular the well-known following one.
\begin{theo}[$\aleph$-dcpo completion]  
Given a preorder $A$, $\adcpo(A)$ and the order preserving 
map $i_A: A \rightarrow \adcpo(A)$ sending $a$ to the downset 
generated by $a$, satisfy the following universal property.
Composing with $i_A$ defines an equivalence of preorders:
\begin{center}
$- \circ i_A : [\adcpo(A),B]' \cong [A,B]$
\end{center}
where:
\begin{itemize}
\item $[A,B]$ is the pointwise preorder of monotonous
maps $A \rightarrow B$;
\item $[\adcpo(A),B]'$ is the pointwise 
preorder of maps $\adcpo(A) \rightarrow B$ that  
preserve $\aleph$-directed least upper bounds.
\end{itemize}
\end{theo}
 
The case $\aleph = \omega$ is quite popular in Computer
Science where our $\omega$-dcpos and $\omega$-dcpo morphisms
are called respectively {\em dcpos} and {\em continuous maps}.
%

\begin{thebibliography}{}
%
%
\bibitem{ABLR02}
J. Adamek, F.Borceux, S.Lack, J.Rosicky:
A classification of accessible categories,
{\it J. Pure Appl. Algebra} \textbf{175} (2002) 7-30.
\bibitem{AK88}
M.H. Albert, G.M. Kelly:
The closure of a class of colimits,
{\it J. Pure Appl. Algebra} \textbf{51} (1988) 1-17.
\bibitem{Bet85}
R. Betti:
Cocompleteness over coverings,
{\it J. Austral. Math. Soc.} \textbf{(Series A) 39} (1985) 169-177.
\bibitem{BCSW83}
R. Betti, A. Carboni, R. Street, R. Walters:
Variation through enrichment,
{\it J. Pure Appl. Algebra} \textbf{29} (1983) 109-127.
\bibitem{Bor94-1}
F. Borceux:
{\it Handbook of categorical algebra 1},
Cambridge University Press (1994).
\bibitem{BvBR98} 
M.M. Bonsangue, F. van Breugel, J.J.M.M Rutten:
Generalized metric spaces: completion, topology, and powerdomains
via the Yoneda embedding, {\it Theoretical Computer Science} 
\textbf{193} (1998) 1-51.
\bibitem{FL82}
P. Fletcher and W. Lindgren:
{\it Quasi-uniform spaces},
Lecture Notes in Pure and Applied Maths - Marcel Dekker ed. (1982).
\bibitem{Fla92}
R.C. Flagg:
Completeness in Continuity Spaces,
{\it AMS Conference proceedings} \textbf{vol. 13} (1992) 183-199.
\bibitem{Fla97}
R.C. Flagg: Quantales and Continuity Spaces,
{\it Algebra Universalis} \textbf{vol. 37} (1997) 257-276.
\bibitem{Ko67}
A. Kock: Limit monads in categories,
{\it Aarhus Univ. Mat. Inst.} Preprint \textbf{No. 6} (1967).
\bibitem{Ko95}
A. Kock: Monads for which structures are adjoint to units,
{\it J. Pure and Appl. Algebra} \textbf{104} (1995) 41-59. 
\bibitem{Kel82}
G.M. Kelly:
{\it Basic concepts of enriched category theory},
London Mathematical Society Lecture Note Series 64
- Cambridge University Press (1982).\\
{\it Reprints in Theory and Applications of Categories} 
\textbf{No.10} (2005)
\bibitem{Kel82-2}
G.M. Kelly:
Structures defined by finite limits in the enriched
context, 1,
{\it Cahiers de Top. et Geo. diff.} \textbf{Vol. XXIII-1} (1982).
\bibitem{KS05}
G.M. Kelly and V. Schmitt:
Notes on categories with colimits of some class,
{\it Theory and Applications of categories} \textbf{vol 14} n.17
(2005) 399-423.
\bibitem{KuSc02}
H.P Kuenzi and M.P.Schellekens:
On the Yoneda completion of a quasi-metric space,
{\it Theoretical Computer Science} \textbf{vol 276} (2002).
\bibitem{Law73}
F.W. Lawvere:
Metric spaces, generalized logic, and closed categories,
{\it Rend. Sem. Mat. Fis. di Milano} \textbf{43} (1973) 135-166.\\
{\it Reprints in Theory and Applications of Categories} 
\textbf{No. 1} (2002).
\bibitem{Law86}
F.W. Lawvere: Taking categories seriously,
{\it Revista Colombiana de Matem\`aticas} \textbf{XX} (1986) 147-178.\\
{\it Reprints in Theory and Applications of Categories} 
\textbf{No. 8} (2005).
\bibitem{Rut96}
J.M.M.M. Rutten:
Elements of generalized ultrametric domain theory,
{\it Theoretical Computer Science} \textbf{70} (1996) 349-381.
\bibitem{Sch03}
V. Schmitt: Enriched categories and quasi-uniform spaces,
{\it Electronic Notes in Theoretical Computer Science} 
\textbf{vol.73.} (2004) 165-205.\\
http://www.elsevier.nl/locate/entcs/volume73.html
\bibitem{Str83-1}
R. Street: Absolute colimits in enriched categories,
{\it Cahiers de Top. et Geo. diff.} \textbf{Vol. XXIV-4} (1983).
\bibitem{Str83-2}
R. Street:
Enriched categories and cohomology,
{\it Quaestiones Math.} \textbf{6} (1983), 265-283.
\bibitem{Tho82} 
W. Tholen: 
Completions of categories and shape theory, 
{\it Seminarberichte aus dem Fachbereich
Mathematik und Informatik der Feruniversitat Hagen} \textbf{12},
(1982) 125-142.
\bibitem{Vic}
S. Vickers: Localic completions of quasimetric spaces,
{\it Preprint} (2003).
\bibitem{Vic05}
S. Vickers: Localic completions of quasimetric spaces,
{\it Theory and Applications of Categories} \textbf{14}
(2005) 328-356.
\bibitem{Wag94}
K.R. Wagner:
Solving recursive domain equations with enriched
categories,
PhD Thesis, Carnegie Mellon University, Pittsburgh, July 1994.
{\it Technical report CMU-CS-94-159}.
\bibitem{Wal81}
R.F.C Walters:
Sheaves and Cauchy-complete categories,
{\it Cahiers de Topo. et G\'eom. Diff. Cat\'egoriques}
\textbf{13} (1981) 217-264. 
\bibitem{Wal82}
R.F.C Walters:
Sheaves on sites as Cauchy-complete categories,
{\it J. Pure and Appl. Algebra} \textbf{24} (1982) 95-102. 
\bibitem{Wo78}
R.J. Wood: Free colimits,
{\it J. Pure and Appl. Algebra} \textbf{10} (1978) 73-80. 

\end{thebibliography}
%

\end{document}